\documentclass{amsart}

\newtheorem{theorem}{Theorem}
\newtheorem{lemma}[theorem]{Lemma}

\newtheorem{definition}[theorem]{Definition}
\newtheorem{example}[theorem]{Example}
\newtheorem{remark}[theorem]{Remark}

\def\QED{\quad\blackslug\lower 8.5pt\null}

\begin{document}

\begin{center}
{\Large \bf An affine analogue of the

\vspace*{1mm}

Hartman-Nirenberg cylinder theorem}

\vspace*{3mm}

 {\large
 Maks A. Akivis and  Vladislav V. Goldberg}

\end{center}
\vspace*{3mm}

{\footnotesize
\textbf{Abstract.} Let $X$ be a smooth, complete, connected
 submanifold of dimension $n < N$ in a complex affine space
$A^N ({\mathbb C})$, and $r$ is the rank of its Gauss map $\gamma, \,
\gamma (x) = T_x (X)$.  The authors prove that if $2 \leq r \leq n - 1,
\, N - n \geq 2,$
and in the pencil of the second fundamental forms
of $X$, there are two forms defining a regular
pencil all eigenvalues of which are distinct, then the submanifold
$X$ is a cylinder with  $(n-r)$-dimensional
plane generators erected over
a smooth, complete, connected  submanifold $Y$ of rank $r$
and dimension $r$. This result is an affine analogue
of the Hartman-Nirenberg cylinder
theorem proved for $X \subset R^{n+1}$ and $r = 1$.
For $n \geq 4$ and $r = n - 1$, there exist complete connected
 submanifolds $X \subset A^N ({\mathbb C})$ that are not cylinders.

\vspace*{5mm}


\vspace*{3mm}


}

\setcounter{equation}{0}

\vspace*{5mm}

\textbf{1. The Hartman-Nirenberg cylinder theorem.}
The Hartman-Nirenberg cylinder theorem in  an $(n+1)$-dimensional Euclidean
space $R^{n+1}$ was first proved  in  \cite{HN}. This theorem
states:

\noindent
\textit{Let $S: M \rightarrow R^{n+1}$ be  a connected, $C^2$, orientable
hypersurface in an $(n+1)$-dimensional space $R^{n+1}$.
 If $S$ is of  constant  zero curvature, then it is an $(n-1)$-cylinder
 $($i.e., an $n$-dimensional cylinder with $(n-1)$-dimensional generators
 erected over a curve$)$ in the
sense that $S$ has a parametrization $($in the large$)$ of the form}
\begin{equation}\label{1}
\vec{v} = \vec{v} (x) = \sum_{i = 1}^{n-1} \vec{a}_j x^j + \vec{b} (x^n)
\; \mbox{{\em for all}} \;\; x = (x^1, \ldots, x^n),
\end{equation}
\noindent
{\em where $\vec{a}_1, \ldots, \vec{a}_{n-1}$ are constant
vectors in $R^{n+1}$, $\vec{b} (x^n)$ is a vector-valued
function of a variable $x^n$ of class $C^2$ in $R^{n+1}$, and
$\vec{a}_1, \ldots, \vec{a}_{n-1}, \,\partial \vec{b}/\partial
x^n$ is a set of orthonormal vectors.}

In the proof of this theorem, the authors first proved
that the vanishing of the Gaussian curvature
implies that  the rank $r (x)$ of the Gaussian map of $S$ does not exceed
one, $r (x) \leq 1$. If  $r (x) = 0$, then $S$ is a hyperplane.
In the case  $r (x) = 1$,
$S$ is an $(n-1)$-cylinder
which can be parametrized as indicated in equation (1).

The proof of this theorem in \cite{HN} is based on
the lemma on the constancy of a certain unique $(n-1)$-plane
from the paper \cite{CL} of Chern and Lashof.
Sternberg \cite{Ste} called this lemma
 Lemma of Chern-Lashof-Hartman-Nirenberg.
A projective analogue of this lemma is in
 \cite{AG1}, Theorem 4.1 (see also \cite{AG2}, Theorem 1).

Note that in \cite{HN} and \cite{Ste}, the authors obtain
an $(n-1)$-cylinder, i.e., a cylinder in $R^{n+1}$
with $(n-1)$-dimensional plane generators erected
over a curve. The reason that they did not get an
$(n-r)$-cylinder, i.e., an $n$-dimensional cylinder in $R^{n+1}$
with $(n-r)$-dimensional plane generators erected
over an $r$-dimensional submanifold, where $r = 1, \ldots , n-1$,
is that their theorem hypotheses imply that the
rank $r (x)$ of the Gaussian map of $S$ does nor exceed one.

The authors noted in  \cite{HN}: ``A similar result under weaker
differentiability hypotheses  has been   stated
for $n = 2$ by Pogorelov'' (see  \cite{P1} and \cite{P2}).
For $R^3$, the results of  \cite{HN} were developed further by Stocker
(see    \cite{Sto1} and \cite{Sto2}).

In recent papers (see for example,  \cite{I1} and \cite{I2})
the authors state the Hartman-Nirenberg cylinder theorem
by saying that ``{\em a properly embedded developable
hypersurface in $R^{n+1}$ of $\mbox{{\rm rank}} \, (\gamma)
\leq 1$ is necessarily a cylinder}''.
A similar result is known for a complex Euclidean space $C^{n+1}$
(see \cite{Ab}; see also the survey \cite{B}).

The authors would like to thank Professor Go-o Ishikawa
for very useful discussions in the process of writing
this paper.
We also would like to thank J. Piontkowski who indicated an
inaccuracy in the first version of our paper.

\textbf{2. An affine analogue of the Hartman-Nirenberg cylinder theorem.}
The Hartman-Nirenberg cylinder theorem is of affine nature. In fact,
the notion of a cylinder appearing in the theorem
conclusion is an affine notion. As to the theorem hypotheses,
although the notion of the Gaussian curvature is not affine, the
notion of the rank $r$ of the Gauss map, which is fundamental in
the proof of this theorem and whose boundedness, $r (x) \leq 1,$
is implied by the vanishing of the Gaussian curvature,
is even of projective nature. This is a reason that it is
interesting to consider an affine analogue of the
Hartman-Nirenberg cylinder theorem.

We recall that an $l$-cylinder $X$ in an  affine space $A^N$
over the field of complex or real numbers is
defined as a smooth $n$-dimensional submanifold bearing
$l$-dimensional plane generators, $l < n$, which are
parallel one to another. An $l$-cylinder is a tangentially
degenerate submanifold of rank $r = n - l$.
In an affine space $A^N,  N > n$,
   an $l$-cylinder can be defined by a parametric equation
\begin{equation}\label{2}
  \vec{\phi} (y^1, \ldots , y^n) = \sum_{j = 1}^l \vec{a}_j y^j +
 \vec{b} (y^{l+1}, \ldots , y^n),
\end{equation}
where $\vec{a}_j$ are constant vectors in $A^N$
and $\vec{b} (y^{l+1}, \ldots , y^{n})$
is a vector-valued function of $r = n - l$ variables
defining in $A^N$ a director variety $Y$ of the cylinder
$X$, and the vectors $\vec{a}_j$ and
$\vec{b}_p = \frac{\partial \vec{b}}{\partial y^p}$ are
linearly independent.

The aim of this paper is to prove the following \textit{affine
analogue of the Hartman-Nirenberg cylinder theorem}.

\begin{theorem}
Let $X^n$ be a smooth, complete, connected  submanifold
 of constant rank $r, \; 2 \leq r \leq n - 1,$
in a complex affine space $A^N (\mathbb{C}), \, N - n \geq 2$.
Suppose that in the pencil of the second fundamental forms
of $X$, there are two forms defining a regular
pencil all eigenvalues of which are distinct.
 Then the submanifold $X$  is a cylinder with $l$-dimensional
 plane generators, $l = n - r \geq 2$,
 and an $r$-dimensional tangentially nondegenerate
  director variety  $Y$. In $A^N (\mathbb{C})$ such a cylinder
 can be defined by parametric equation $(2)$.
\end{theorem}

We recall that the eigenvalues of a pencil of quadratic forms
defined by the forms
$$
\phi' =  b'_{pq} \theta^p \theta^q \; \text{and} \;
\phi'' =  b''_{pq} \theta^p \theta^q
$$
are determined by the characteristic equation
$$
 \det (b''_{pq} + \lambda b'_{pq}) = 0.
$$

Note also that if a submanifold $X \subset A^N (\mathbb{C})$
is smooth (i.e., differentiable), then it is complex-analytic, and
in particular, it is twice differentiable.

Note that the condition $N - n \geq 2$ in our cylinder theorem
 is essential since if $N - n = 1$, in the space $A^{n+1}
 (\mathbb{C})$,  there exist complete, connected, complex-analytic
 hypersurfaces $X^n$ of rank $r < n$ without singularities that
 are not cylinders. An example of such hypersurface
 $X^3$ in the space $A^4 (\mathbb{C})$ was considered
 by Sacksteder in \cite{S}. An algebraic ruled noncylindrical hypersurface
without singularities in $A^4 (\mathbb{C})$
 was constructed by Bourgain and published by Wu
 \cite{W} (see also \cite{I1} and \cite{I2}). This example
 was studied in detail in \cite{AG3}.


Note that an affine cylinder theorem in other formulations
was presented in the paper \cite{NP} by Nomizu and Pinkall
(see also the book \cite{NS}) and in the papers
\cite{O1}, \cite{O2}, and \cite{O3} by Opozda.
Their  affine cylinder theorems
 give sufficient conditions for a hypersurface
(i.e., a submanifold of codimension one)
$X$ in $A^{n+1}$ to be a cylinder erected over a curve
with $(n-1)$-dimensional plane generators.
Our  affine cylinder theorem gives  sufficient conditions for a
submanifold $X$ of any codimension and any rank
$r, \, 2 \leq r \leq n - 1$, in $A^N, \, N - n \geq 2, $
to be a cylinder erected over a submanifold
of dimension $r$ and rank $r$
with $(n-r)$-dimensional plane generators.
In a recent paper \cite{Pi}, Piontkowski considered in $P^N$ complete
varieties with degenerate Gauss maps with rank equal 2 or 3 or 4
and with all singularities
located at a hyperplane at infinity. In particular, as an extreme
case, he obtained an affine cylinder theorem
for varieties of rank 1 and any codimension.
So, our affine cylinder theorem
for submanifolds of codimension greater than 2 and
rank $r \geq 2$ complements substantially all previously
known  affine cylinder theorems which were for hypersurfaces
of rank 1.

Theorem 1 is a simple consequence of a certain
structure theorem (Lemma 2) which is of projective nature
(see \cite{AG2}, where one can also find other
structure theorems for submanifolds with
degenerate Gauss map). In order to present this theorem,
we need to introduce some notions and notations.

\textbf{3. Basic equations and focal images.}
An $n$-dimensional submanifold $X$ of a projective space $P^N
(\mathbb{C})$ is called  {\em tangentially degenerate} if the rank of its
Gauss mapping $\gamma: X \rightarrow G (n, N)$ is less than $n, \;
r = \mbox{{\rm rank}} \; \gamma < n$.
Here $x \in X, \; \gamma (x) = T_x (X)$, and $T_x (X)$ is the
tangent subspace to $X$ at $x$ considered as an $n$-dimensional
projective space $P^n (\mathbb{C})$.  The number  $r$
is also called the {\em rank} of $X, \; r  =  \mbox{{\rm rank}} \; X$.
We assume that the rank $r$ is constant on the submanifold $X$.

Denote by $L$ a leaf of this map, $L =
 \gamma^{-1} (T_x) \subset X$.
It is easy to prove that a leave $L$ of the Gauss map
$\gamma$ is a subspace of $P^N (\mathbb{C}), \; \dim L = n - r = l$
or its open part (see for example, \cite{AG1}, p. 115,  Theorem 4.1).

The foliation on $X$ with leaves $L$ is called
the {\em Monge--Amp\`{e}re foliation}.
As we did in \cite{AG2}, in this paper we extend the leaves
of the Monge--Amp\`{e}re foliation to a projective space $P^l ({\mathbb C})$
assuming that $L \sim P^l ({\mathbb C})$ is a plane generator of the
submanifold $X$. In this sense
the submanifold $X$ is complete.
 As a result, we have $X = f (P^l ({\mathbb C})\times M^r)$,
where $M^r$ is a complex  variety of parameters,
and $f$ is a differentiable map $f:
P^l ({\mathbb C}) \times M^r \rightarrow P^N ({\mathbb C})$.
The Monge--Amp\`{e}re foliation
is locally trivial but, as we will see later,
its leaves $L$ can have singularities.

In $P^N (\mathbb{C})$, we consider a manifold of projective frames
$\{A_0, A_1, \ldots , A_N\}$. On this manifold
\begin{equation}\label{eq:3}
 d A_u =  \omega_u^v A_v, \;\; u, v, = 0, 1, \ldots , N,
\end{equation}
where the sum $\omega_u^u = 0$. The 1-forms
$\omega_u^v$ are linearly expressed in terms of
the differentials of parameters of the group
of projective transformations of the space $P^N (\mathbb{C})$.
These 1-forms
satisfy the structure equations
\begin{equation}\label{eq:4}
 d  \omega_u^v  =  \omega_u^w \wedge  \omega_w^v
\end{equation}
of the space $P^N (\mathbb{C})$ (see, for example, \cite{AG1}, p. 19).
Equations (4) are the conditions of complete integrability
of equations (3).

Consider a  tangentially degenerate submanifold $X \subset P^N (\mathbb{C}),
\; \dim X = n, \; \mbox{{\rm rank}}\; X = r < n$. In addition,
as  above, let $L$ be a plane generator of the manifold
$X, \; \dim L = l = n - r$; let $T_L, \; \dim T_L = n,$ be the tangent
subspace to $X$ along the generator $L$, and let $M$ be a base
manifold for $X, \; \dim M = r$. Denote by $\theta^p, \; p = l + 1,
\ldots, n$,  basis forms on the variety $M$. These forms satisfy
the structure equations
\begin{equation}\label{eq:5}
 d  \theta^p  =  \theta^q \wedge  \theta_q^p, \;\; p, q = l + 1,
 \ldots , n,
\end{equation}
of the variety $M$. Here $\theta_q^p$ are 1-forms defining
transformations of first-order frames on $M$.

For a point $x \in L$, we have $d x \in T_L$. With $X$, we
associate a bundle of projective frames $\{A_i, A_p, A_\alpha\}$
such that $A_i \in L, \; i = 0, 1, \ldots, l; \; A_p \in T_L, \;
p = l + 1, \ldots , n$. Then
\begin{equation}\label{eq:6}
\renewcommand{\arraystretch}{1.3}
\left\{
\begin{array}{ll}
dA_i =  \omega_i^j A_j + \omega_i^p A_p, \\
dA_p =  \omega_p^i A_i + \omega_p^q A_q +
 \omega_p^\alpha A_\alpha, \; \alpha = n + 1, \ldots, N.
\end{array}
\right.
\renewcommand{\arraystretch}{1}
\end{equation}
It follows from the first equation of (6) that
\begin{equation}\label{eq:7}
  \omega_i^\alpha = 0.
\end{equation}
Since for $\theta^p = 0$ the subspaces $L$ and $T_L$ are fixed, we
have
\begin{equation}\label{eq:8}
  \omega_i^p = c^p_{qi} \theta^q, \;\;   \omega_p^\alpha
= b^\alpha_{qp} \theta^q.
\end{equation}

Since the manifold of leaves $L \subset X$ depends on $r$
essential parameters, the rank of the system of 1-forms
$\omega^p_i$ is equal to $r$, $\mbox{{\rm rank}} \; (\omega^p_i) =
r$. Similarly, we have $\mbox{{\rm rank}} \; (\omega^\alpha_p) =
r$.

 Denote by
$C_i$ and $B^\alpha$ the  $r \times r$ matrices occurring in
equations (8):
$$
C_i = (c^p_{qi}), \;\; B^\alpha = (b_{qp}^\alpha).
$$
These matrices are defined in a second-order neighborhood
of the submanifold $X$.

Exterior differentiation of equations (7) by
means of  structure equations (4) leads to the exterior
quadratic equations
$$
\omega_p^\alpha \wedge \omega_i^p = 0.
$$
Substituting expansions (8) into the last equations, we find that
\begin{equation}\label{eq:9}
b^\alpha_{qs} c^s_{pi} = b^\alpha_{ps} c^s_{qi}.
\end{equation}
Equations (7), (8), and (9) are called the {\em basic equations}
in the theory of tangentially degenerate submanifolds.

Relations (9) can be written in the matrix form $$ (B^\alpha
C_i)^T = (B^\alpha C_i), $$ i.e., the matrices $$ H^\alpha_i =
B^\alpha C_i = (b^\alpha_{qs} c^s_{pi}) $$ are symmetric.

Let $x = x^i A_i$ be an arbitrary point of a leaf $L$. For such a
point we have
$$
d x = (dx^i + x^j \omega_j^i) A_i + x^i
\omega_i^p A_p.
$$
It follows that
$$
d x \equiv (A_p c^p_{qi}
x^i) \theta^q \pmod{L}.
$$
The tangent subspace $T_x$ to the
manifold $X$ at a point $x$ is defined by the points $A_i$ and $$
\widetilde{A}_q (x) = A_p c^p_{qi} x^i, $$ and therefore $T_x
\subset T_L$.

A point $x$ is a {\em regular} point of a leaf $L$ if
$T_x = T_L$. The regular points are determined by the condition
\begin{equation}\label{eq:10}
J (x) = \det (c^q_{pi} x^i) \neq 0.
\end{equation}
If $J (x) = 0$ at a point $x$, then $T_x$ is a proper subspace
of $T_L$, and a point $x$ is said to be a {\em singular point}
of a leaf $L$.

The determinant (10) is the Jacobian of the map $f: P^l \times
M^r \rightarrow P^N$.
The singular points of a leaf $L$ are determined by the
condition $J (x) = 0$. In a leaf $L$, they form an algebraic
submanifold of dimension $l - 1$ and degree $r$. This
hypersurface (in $L$) is called the {\em focus hypersurface}
and is denoted by $F_L$.
By (10), the equations $J (x) = 0$ of the focus hypersurface
in the plane generator $L$ of the manifold $X$ can be written as
\begin{equation}\label{eq:11}
 \det (c^q_{pi} x^i) = 0.
\end{equation}

For the point $x = x^i A_i \in L$, by (6) we have
$$
d^2 x \equiv A_\alpha b^\alpha_{qs} c^s_{pi} x^i  \theta^p \theta^q
\pmod{T_L}.
$$
Thus, the points
$$
A_{pq} = A_\alpha b^\alpha_{qs} c^s_{pi} x^i, \;\;
A_{pq} = A_{qp},
$$
together with the points $A_i$ and $A_p$ define the
osculating subspace $T^2_x (X)$. Its dimension is
$$
\dim T_x^2 (X) = n + m,
$$
where $m$ is the number of linearly independent points
among the points
$A_{pq}, \; m \leq \mbox{{\rm min}} \{\frac{r (r + 1)}{2}, N - n\}$.
But since at a regular point $x \in L$  condition (10) holds,
the number $m$ is the number of linearly independent points
among the points
$$
\widetilde{A}_{pq} = A_\alpha b^\alpha_{pq}.
$$
Hence,  the osculating subspace $T^2_L (X)$
is the same for all regular points $x \in L$, and we
denote it by $T_L^2 (X)$.

We further assume that the point $A_0$ is a regular point
of the generator $L$. Then
equation (10) implies that  $\det (c_{q0}^p) \neq 0$, and
the 1-forms $ \omega_o^p = c^p_{q0} \theta^q$
are linearly independent.  Thus, they
can be taken as basis forms of the manifold $M^r$, i.e.,
we can set $\omega_0^p = \theta^p$. Then
the matrix $C_0 = (c_{q0}^p)$ becomes the
identity matrix, $c_{q0}^p = \delta_q^p$.
This and equations (9) imply that
$$
b^\alpha_{pq} = b^\alpha_{qp},
$$
i.e., the matrices $B^\alpha$ become symmetric.

The number $m$ coincides with the number of linearly independent
matrices $B^\alpha = (b^\alpha_{pq} )$ and with the number of linearly independent
second fundamental forms
$$
\phi^\alpha = b^\alpha_{pq}  \omega_0^p  \omega_0^q.
$$
of the submanifold $X$ at its regular point $A_0$.

\textbf{4. Decomposition of the focus hypersurface $\boldsymbol{F_L}$.}
We prove now the following lemma;
\setcounter{theorem}{0}

\begin{lemma}
Suppose that  $m \geq 2$, and  in the pencil of the second fundamental forms
$\xi_\alpha \phi^\alpha$ of $X$, there are two forms $\phi'$ and
$\phi''$ defining a regular pencil, and  all eigenvalues
$\lambda_p, p = l + 1, \ldots, n,$ of this pencil
defined by the characteristic equation
\begin{equation}\label{eq:12}
 \det (b''_{pq} + \lambda b'_{pq}) = 0
\end{equation}
are  distinct. Then all the matrices $C_i = (c^p_{qi})$
can be diagonalized  simultaneously with the matrices
 $B' = (b'_{pq})$ and $B'' = (b''_{pq})$.
\end{lemma}

{\sf Proof.} In fact, by the lemma condition, the matrices
 $B'$ and $B''$ can be simultaneously   diagonalized:
\begin{equation}\label{eq:13}
 B' = (b'_{pp}), \;\; B'' = (b''_{pp}).
\end{equation}
Since all roots of equation (12) are distinct, then
\begin{equation}\label{eq:14}
\frac{b'_{pp}}{b''_{pp}} \neq \frac{b'_{qq}}{b''_{qq}}, \;\; p \neq q.
\end{equation}

Consider now equations (9). By (13), equations (9) take the form
\begin{equation}\label{eq:15}
b'_{qq} c^q_{pa} = b'_{pp} c^p_{qa}, \;\;
b''_{qq} c^q_{pa} = b''_{pp} c^p_{qa},
\end{equation}
where $a = 1, \ldots, l.$ But by the inequality (14),
equations (15) imply that
\begin{equation}\label{eq:16}
c^p_{qa} = 0, \;\;\; p \neq q.
\end{equation}
Since in addition, we have $C_0 = (\delta_q^p)$, each of the matrices
$C_i$ has a diagonal form. 

Lemma 1 implies the following corollary.

\textbf{Corollary.} \textit{Under the conditions of Lemma $1$, the focus
hypersurface $F_L$ of its generator $L$ decomposes into $r$
hyperplanes defined by the equations
\begin{equation}\label{eq:17}
x^0 + c^p_{pa} x^a  = 0, \;\; a = 1, \ldots, l,
\end{equation}
where there is no summation over $p$.}

 In fact, by Lemma 1, the matrices $C_0$ and $C_a$
can be simultaneously  diagonalized, and
$C_0 = (\delta_q^p)$. As a result, equation (11) of the
focus hypersurface $F_L$ takes the form
\begin{equation}\label{eq:18}
 \prod_{p = l+1}^n  (x^0 + c^p_{pi} x^i) = 0.
\end{equation}
Therefore, the hypersurface $F_L$  decomposes into $r$
hyperplanes (17).

We can reformulate Lemma 1 in geometric terms by using
the notion of the focus hypercone $\Phi_L$ formed
by singular tangent hyperplanes of the submanifold $X$
containing the tangent subspace $T_L$. The hypercone $\Phi_L$
is determined by the equations
$$
\det (\xi_\alpha b^\alpha_{pq}) = 0
$$
(see \cite{AG1}, p. 119, or \cite{AG2}). The new
formulation of Lemma 1 is:

\emph{If $m \geq 2$, and the focus hypercone $\Phi_L$
does not have multiple components, then
 the focus hypersurface $F_L$ of  its generator $L$ decomposes into $r$
hyperplanes $F_p$ lying in $L$.}

\textbf{5. Sufficient conditions for $\boldsymbol{X}$ to be a cone.}
Lemma 1 allows us to prove the following lemma
giving sufficient conditions for a submanifold $X$
with a degenerate Gauss map to be a cone.

\begin{lemma}
Suppose that  a tangentially degenerate submanifold $X$
satisfies the following conditions:
\begin{itemize}
\item[(i)] The conditions of Lemma 1.
\item[(ii)] $r \geq 2$.
\item[(iii)] All hyperplanes $F_p$, into which by Lemma $1$
the focus hypersurface $F_L$ decomposes, coincide.
\end{itemize}
Then the  submanifold $X$ is a cone with $(l-1)$-dimensional
vertex and $l$-dimensional plane generators.
\end{lemma}

{\sf Proof.} In fact, if all hyperplanes $F_p$ defined by
equations (17) coincide, then by placing the points $A_b, b = 1,
\ldots , l$, of our moving frame
into this hyperplane, we reduce the equation of this hyperplane to
the form $x^0 = 0$. As a result, all the coefficients $c^p_{pi}$
vanish. Since we also have conditions (16), it follows that
$$
c^p_{qa} = 0, \;\; a = 1, \ldots , l,
$$
for all $p, q = l + 1, \ldots , n$.
Now equations (8) imply that
\begin{equation}\label{eq:19}
 \omega_a^p  = 0.
\end{equation}
Taking exterior derivatives of equations (19), we obtain
$$
\omega_a^0 \wedge \omega_0^p  = 0.
$$
Since the forms $\omega_0^p$ are linearly independent (they are
basis forms of the variety $M^r$) and $r \geq 2$,
it follows that
\begin{equation}\label{eq:20}
 \omega_a^0  = 0.
\end{equation}

Equations (18) and (19) imply that
$$
d A_b =  \omega_b^c A_c,
$$
i.e., the $(l-1)$-plane $A_1 \wedge \ldots \wedge A_l$ is
constant. Thus, the submanifold $X$ is a cone with $(l-1)$-dimensional
vertex.

\textbf{6. Proof of Theorem 1.}
We have a submanifold  $X^n$ of dimension $n$
 and rank $r,\; \mbox{{\rm rank}} \; X = r,$ without
singularities in an  affine space $A^N (\mathbb{C}), \,
N - n \geq 2$. Thus all
 singular points of $X$ are located at
 a hyperplane at infinity  $P^{N-1}_\infty \subset
  \overline{A}^N = A^N \cup
 P^{N-1}_\infty$.

 Now for $r \geq 2$, we can apply Lemma 2.
The cones in $P^N (\mathbb{C})$  are projectively equivalent to
the cylinders $X^n \subset A^N ({\mathbb C})$:
a cylinder $X^n \subset  A^N ({\mathbb C})$
is obtained from a cone $\widetilde{X}^n \subset P^N (\mathbb{C})$
by means of a projective transformation sending
the cone vertex $F, \, \dim \, F  = l - 1$,
into $P^{N-1}_\infty \subset A^N ({\mathbb C})$.
The director  variety $Y, \; \dim Y
= \mbox{{\rm rank}} \; Y = r$, of the cylinder $X^n$
lies in the proper subspace $T, \, \dim \, T = N - l;
 T$  is complementary to $F$. Thus, by Lemma 2,
 $X^n$ is an $(n-r)$-cylinder  erected over an $r$-dimensional
 submanifold (a director variety) $Y$
 and having $(n-r)$-dimensional plane generators.

\makeatletter \renewcommand{\@biblabel}[1]{\hfill#1.}\makeatother
\bibliographystyle{amsplain}


\noindent {\em Authors' addresses}:\\

\noindent
\begin{tabular}{ll}
M.~A. Akivis &V.~V. Goldberg\\ Department of Mathematics
&Department of Mathematical Sciences\\ Jerusalem College of
Technology---Mahon Lev &  New Jersey Institute of Technology \\
Havaad Haleumi St., P. O. B. 16031 & University Heights \\
 Jerusalem 91160, Israel &  Newark, N.J. 07102, U.S.A. \\
 & \\
 E-mail address: akivis@avoda.jct.ac.il & E-mail address:
 vlgold@m.njit.edu
 \end{tabular}

\end{document}